# Phasor-Based Approach for Harmonic Assessment from Multiple Distributed Energy Resources

Reza Arghandeh, Alexandra von Meier
Member, IEEE
CIEE, EECS Dept
University of California, Berkeley
Berkeley, CA, USA
arghandeh@berkeley.edu

Robert Broadwater
Senior Member, IEEE
ECE Dept
Virginia Tech
Blacksburg, VA, USA

*Abstract*— **This paper analyzes impacts and interactions of harmonics from multiple sources, especially distributed energy resources, on distribution networks. We propose a new index, the Phasor Harmonic Index (PHI), that considers both harmonic source magnitude and phase angle, while other commonly used harmonic indices are based solely on magnitude of waveforms. The use of such an index becomes feasible and practical through emerging monitoring technologies like micro-Synchrophasors in distribution networks that help measure and visualize voltage and current phase angles along with their magnitudes. A very detailed model of a distribution network is also needed for the harmonic interaction assessment in this paper.**

*Nomenclatures—*

| | |
|---|---|
| $V_h$ | : Voltage magnitude for frequency order h |
| $I_h$ | : Current magnitude for frequency order h |
| $V_{Total}$ | : Total value for voltage magnitude |
| $I_{Total}$ | : Total value for current magnitude |
| $V_h^{Ph}$ | : Phasor form of voltage for frequency order h |
| $I_h^{Ph}$ | : Phasor form of current for frequency order h |
| $\theta h$ | : Harmonic current phase angle |
| $\varphi_h$ | : Harmonic voltage phase angle |
| $h$ | : Harmonic frequency order |
| $I_{hP}$ | : In-phase current for frequency order h |
| $I_{hQ}$ | : In-quadrature current for frequency order h |
| $V_{hP}$ | : In-phase voltage for frequency order h |
| $V_{hQ}$ | : In-quadrature voltage for frequency order h |
| THDV | : Voltage Total Harmonic Distortion |
| THDI | : Current Total Harmonic Distortion |
| PHI-V | : Phasor Harmonics Index- Voltage |
| PHI-I | : Phasor Harmonics Index - Current |

## I. INTRODUCTION

Harmonics originating from increasingly prevalent nonlinear, power electronic-based loads as well as inverters associated with distributed energy resources (DER) create concerns for distribution network operators. Such harmonics (i.e., multiples of the 50- or 60-Hz fundamental a.c. frequency) can propagate and distort voltage and current waveforms in different parts of distribution networks, resulting in undesirable heating and energy losses (and potentially equipment malfunctions). Harmonics generated by different sources can also interact to either increase or decrease the effects of harmonics, not unlike the familiar constructive or destructive interference of propagating waves[1].

In general, the harmonic impact on power systems is a well-researched topic. Harmonic measurement and filtering in power systems are discussed in [2]. However, existing research mostly considers harmonics as a local phenomenon with local effects [3]. Some papers focus on DER harmonics [4, 5], and some specifically on harmonic filter design for DER units [6, 7]. However, the proposed solutions are local approaches for controlling each inverter. Less literature investigates the impact of harmonic propagation in distribution networks. Harmonic distortion in different distribution transformer types is analyzed in [8], and [9] conducts a sensitivity analysis to find vulnerable buses in distribution networks. However, the authors use the Thévenin equivalent model at each bus instead of the full topological model of the circuit. The impact of aggregated harmonics from DER units in distribution networks is shown in [10], but using single-phase equivalent line models without considering multi-phase line models.

This paper investigates a phasor-based harmonic assessment method to understand the interactions of multiple harmonic sources. The analysis benefits from the detailed imbalanced and asymmetrical distribution network model employed. This model has large numbers of single phase, and multi-phase loads for more realistic harmonic propagation simulations, a level of detail not addressed in the previous harmonic analysis literature. Thus, we can effectively study how DER inverters or other sources may interact to either decrease or increase harmonic distortion throughout the distribution network.

In terms of harmonic distortion quantization, the Total Harmonic Distortion (THD) is the most common index in standards and literature [11, 12]. However, THD is based only on the magnitude of the distorted waveforms. In this paper a new index is proposed called the Phasor Harmonic Index (PHI). The PHI incorporates both magnitude and phase angle information in evaluating distorted waveforms resulting from the phasor-based interaction of multiple harmonic sources. The phasor-based harmonic assessment in this paper is inspired by the micro-Synchrophasor (μPMU) measurement technology in whose development two of the authors are involved [13, 14]. The μPMU is the GPS enabled time-synchronized phasor measurement devices that can measure voltage and current signals till 50th harmonic order.

The paper is organized as follows: Section 2 discusses harmonics assessment framework; Section 3 presents simulations and results, and Section 4 offers concluding remarks.

## II. Harmonic Assessment Framework

### A. Integrated System Model for Distribution Networks

An Integrated System Model (ISM) is used for distribution network harmonic analysis. Geographical information, component characteristics, load measurements, and supply measurements are included in the model [15, 16]. The ISM model offers a graph-based, edge-edge topology iterator framework that facilitates fast computation times for power flow and other calculations on the large scale networks [17]. The [18] provides further explanation about ISM modeling.

This analysis employs an actual circuit model, as shown in Figure 1. The circuit is 13.2 kV with 329 residential and commercial customers. The model contains unbalanced, single phase and multi-phase loads, and includes distribution transformers and secondary distribution. The two harmonic sources studied are indicated with triangular symbols. The harmonic calculations are presented at two points indicated by arrows in Figure 1. The first point is the substation, and the second point is at the secondary of a distribution transformer located between the two harmonic sources.

### B. Metrics for Vectorial Harmonics Assessments

The summation of harmonic components results in distorted current and voltage waveforms. The most common index used for measuring harmonics in standards and literature is Total Harmonic Distortion (THD) [12]. THD includes the contribution of the magnitude of each harmonic component as given by

$$THDI = \frac{1}{I_1}\sqrt{\sum_{h=2}^{\infty} I_h^2} \times 100 \quad (1)$$

$$THDV = \frac{1}{V_1}\sqrt{\sum_{h=2}^{\infty} V_h^2} \times 100 \quad (2)$$

where *THDI* and *THDV* are THD values for current and voltage, respectively. $I_1$ and $V_1$ are the current and voltage rms values for the fundamental frequency. The definition of THD is addressed in IEEE-519, IEEE-1547, IEC-610000, and EN50160 and other standards for power quality [19, 20].

In some standards, the conventional definition of power factor is modified to account for the contribution of higher frequencies [12], which tend to contribute to total current and thus apparent power, but not to the net transfer of average (real) power. This modified power factor, which no longer assumes a sinusoidal waveform, is called Total Power factor (TPF).

Equation (3) shows the relationship between TPF and THD [21]:

$$TPF = \frac{cos(\delta_1)}{\sqrt{1 + THDI^2}} \quad (3)$$

where $\delta_1$ is the angle between voltage and current at the fundamental frequency. Note that the THD indices are based solely on the magnitude of harmonic components, and the only phase angle difference considered is that between the fundamental voltage and current vectors.

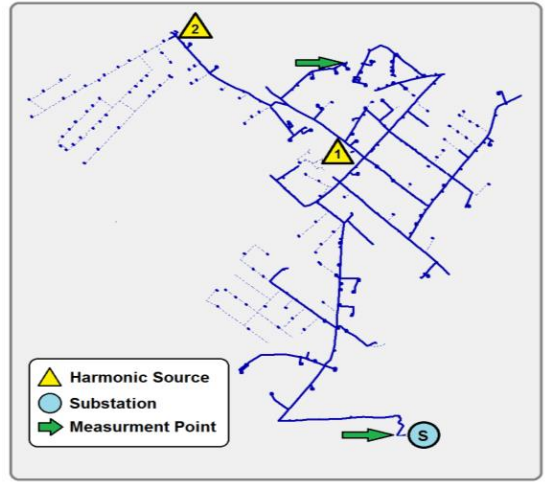

Figure 1. Distrbution Network ISM model for Harmonic Assessment.

Thus the most common indices for harmonic analysis do not account for the phase angles of the harmonic components in harmonic distortion assessment. However, the phase angle of the harmonic waveforms do have an impact on the total distorted current or voltage waveforms. To help quantify the distortion caused by multiple harmonic source interactions, there is a need to consider phase angles of current and voltage waveforms for all frequency orders. The orthogonal form of sinusoidal voltage and current waveforms can be written as:

$$I_{Total} = \sqrt{2}\sum_{h=1}^{\infty} I_h cos(\theta_h)sin(h\omega_0 t) - \sqrt{2}\sum_{h=1}^{\infty} I_h sin(\theta_h)cos(h\omega_0 t) \quad (4)$$

$$V_{Total} = \sqrt{2}\sum_{h=1}^{\infty} V_h cos(\varphi_h)sin(h\omega_0 t) - \sqrt{2}\sum_{h=1}^{\infty} V_h sin(\varphi_h)cos(h\omega_0 t) \quad (5)$$

In (6) and (7), total current and voltage are separated into two in-phase and in-quadrature components. This waveform separation method is similar to some apparent power calculation methods for nonsinusoidal apparent power calculation in [22]. Equations (4) and (5) are rewritten as follows:

$$I_{Total} = \sum_{h=1}^{\infty} I_{hP}\sin(h\omega_0 t) - \sum_{h=1}^{\infty} I_{hQ}\cos(h\omega_0 t) \quad (6)$$

$$V_{Total} = \sum_{h=1}^{\infty} V_{hP}\sin(h\omega_0 t) - \sum_{h=1}^{\infty} V_{hQ}\cos(h\omega_0 t) \quad (7)$$

$$I_{hP} = \sqrt{2}I_h\cos(\theta_h), \qquad I_{hQ} = \sqrt{2}I_h\sin(\theta_h)$$

$$V_{hP} = \sqrt{2}V_h\cos(\varphi_h), \qquad V_{hP} = \sqrt{2}V_h\sin(\varphi_h)$$

The terms $I_{hP}$ and $V_{hP}$ are called in-phase current and voltage components; $I_{hQ}$ and $V_{hQ}$ are called in-quadrature current and voltage components. In some references, in-phase components are called active and in-quadrature components are called non-active components [23].

In this paper, with help of equations (6) and (7), the Phasor Harmonic Index (PHI) is proposed. The PHI is obtained by dividing the summation of in-phase harmonic components by the algebraic sum of harmonic waveform magnitudes as follows:

$$PHI - I = \frac{\sum_{h=1}^{\infty}|I_{hP}|}{\sum_{h=1}^{\infty}|I_h|} \quad (8)$$

$$PHI - V = \frac{\sum_{h=1}^{\infty}|V_{hP}|}{\sum_{h=1}^{\infty}|V_h|} \quad (9)$$

where *PHI-I* and *PHI-V* are Phasor Harmonics Index for current and voltage waveforms, respectively.

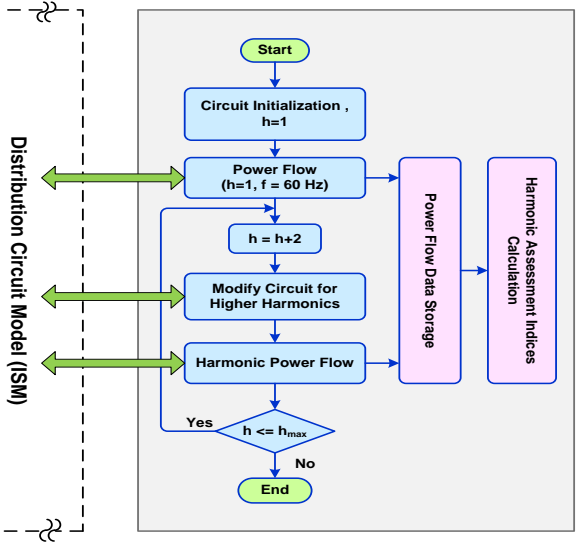

Figure 2. Phasor-Based Harmonic Assessment Diagram

### C. Harmonics Assessment Framework

The harmonic assessment algorithm first uses power flow analysis to calculate the fundamental voltage and current waveforms. Then, the circuit is modified to represent the next higher frequency to be analyzed. The circuit modifications involve revising the network component impedances and the harmonic current injections for the harmonic order to be analyzed. Next the power flow runs to determine the harmonic current and voltage emissions in the circuit for the given harmonic order. This type of analysis is continued until all harmonic orders to be analyzed are completed.

Harmonics are affected by configuration, impedance, and loading of conductors, transformers, and other circuit components[24]. Figure 2 illustrates the flow of the harmonic assessment algorithm. The "Power Flow Data Storage" stores fundamental and higher order harmonic power flow results that are used to calculate harmonic assessment indices. In the work here, $h_{max}$ is 11 and only the odd harmonic orders are taken into account, since these typically dominate.

### III. SIMULATIONS AND HARMONICS ANALYSES

#### A. Assumptions

The research objective is the study of harmonic impacts, apart from whatever technology created the harmonic source. We consider two 3-phase harmonic sources in the distribution network, with equal magnitude on all phases. The harmonic magnitudes are based on test data from actual DER inverters in the field as shown in Table I [15]. The dominant current and voltage harmonic observed through the simulation are of the 3rd, 5th, 7th, 9th and 11th orders. Harmonics of higher orders are neglected due to their small values.

TABLE I. HARMONIC SOURCE MAGNITUDES FROM FIELD TEST DATA [21]

| Frq Order | 3rd | 5th | 7th | 9th | 11th |
|---|---|---|---|---|---|
| Current (p.u %) | 2.83 | 0.52 | 0.84 | 0.21 | 0.03 |

#### B. Phase Angle Impacts on Harmonic Interactions

In systems with multiple harmonic sources, the harmonic distortion interactions are impacted by the vectorial characteristics of the injected harmonic currents. The impact of each harmonic source's phase angle is investigated in this section. For sensitivity analysis purposes, harmonic source phase angles vary as follows: (0°, 15°, 30°, 45°, 60°, 75°, 90°). These angle steps are added to the phase angle sequences for each frequency order [12]. When varying the phase angles of the harmonic sources, the magnitudes of both harmonic sources are maintained as given in Table I.

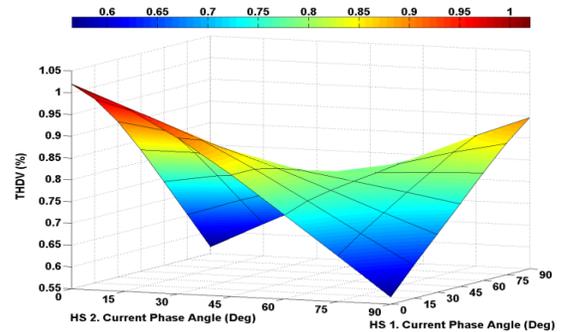

Figure 3. THDV for phase B as a function of harmonic source 1 (HS1) and harmonic source 2 (HS2) phase angles at substation

We begin by considering the THDV at the substation, obtained by vector summation of harmonic contributions from the two sources, on a single phase (say, B). Figure 3 illustrates

THDV at the substation for different combinations of harmonic phase angles at their sources. The THDV surface shapes are similar to the hyperbolic geometrical functions. In this case, the extreme points for three phases occur in 0°, 45°, and 90° phase angles for each harmonic source. Figures 4 shows the analogous THDI surfaces for current in phase B at the substation against phase angles variation over both harmonic sources. There are two near zero points, the harmonic sources cancel out each other and cause the minimum current harmonic distortion, (90°, 0°) and (0°, 90°). This possibility of cancellation is an important observation for multi-source harmonic analysis.

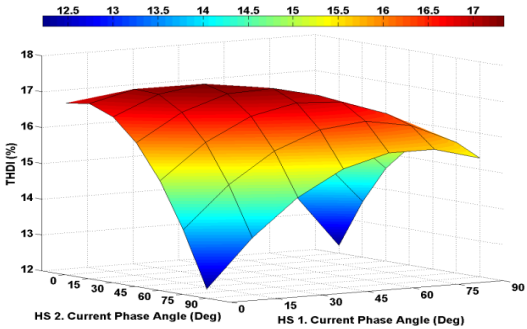

Figure 4. THDI for phase B as a function of harmonic source 1 (HS1) and harmonic source 2 (HS2) phase angles at substation

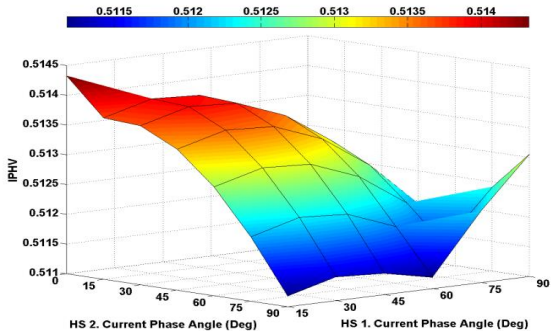

Figure 5. PHI-V for phase B as a function of harmonic source 1 (HS1) and harmonic source 2 (HS2) phase angles at substation

The presented THD sensitivity analysis shows the harmonic sources' critical angles for voltage and current distortion. However, the difference of THD trajectory for voltage and current make it difficult to find the most critical phase angles for harmonic distortion in substations. Figure 5 shows the proposed Index of PHI-V as defined in (9). Because of the contribution of the phase angle in the PHI numerator, PHI contains more information than THD. The PHI values are smaller or equal to 1, because of the "Triangle Inequality" property in a vector space.

### C. Mutual Coupling Impact and Phasor-Based Harmonics

In this section, harmonic sources attached to one phase are analyzed to determine their impact on other phases. The single phase harmonic sources are located at the same place as the three-phase harmonic, and THDV and THDI are again calculated at the substation. With harmonic current injections in only one phase, THDV and THDI indices for the other phases are almost zero. However, mutual couplings do cause distortion in coupled voltage and current waveforms, but these are not reflected in the THDV and THDI indices. The PHI-I index does provide non-zero harmonic distortion values for coupled phases. In Figure 6, the PHI-I values for different phase angles (0° to 90°) of both harmonic sources are classified as a data set presented in the form of a box plot. The box plots show minimum, maximum, mean, and median values of PHI-I calculated for all phases caused by injected harmonics on phase B.

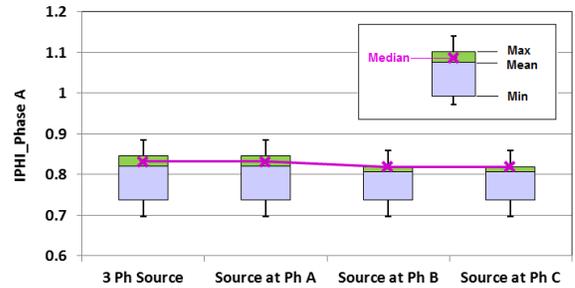

Figure 6. Box plot for PHI-I values with different harmonic source phase angles for the harmonic source injections in only phase B.

These types of sensitivity analyses are not possible with THDV and THDI indices because of the extremely small values of THD in the coupled phases that do not contain the harmonic source.

### D. Phasor-Based Harmonics at Customer Level

To demonstrate the harmonic distortion at different locations on the circuit, harmonic calculations are conducted at a customer load (on the secondary side of a distribution transformer). The customer side measurement point is located between the harmonic sources, as illustrated in Figure 1.

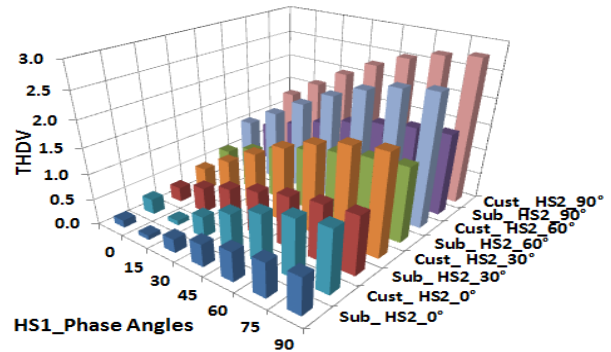

Figure 7. THDV substation vs. customer load point for phase A

Figure 7 compares THDV values at the substation and at the customer load, and Figure 8 presents the THDI values at the substation and at the customer load. It shows the substation experiences more distortion in current than the customer load, while the customer load is exposed to higher harmonic voltage distortion. This result can be understood by considering comparative impedances. Since the impedance looking back into the substation is much smaller than the customer load impedance, a higher portion of harmonic currents flow to the substation than to the customer site. Therefore, the substation has more THDI. The voltage distortion, however, is larger at the customer load. Voltage is the product of impedance and current. The customer load impedance is much larger than the impedance of the path to the substation. The current through the load side is less, but the product of current and impedance for the customer side is higher than for the substation side. Therefore, more voltage distortion is realized at the customer side.

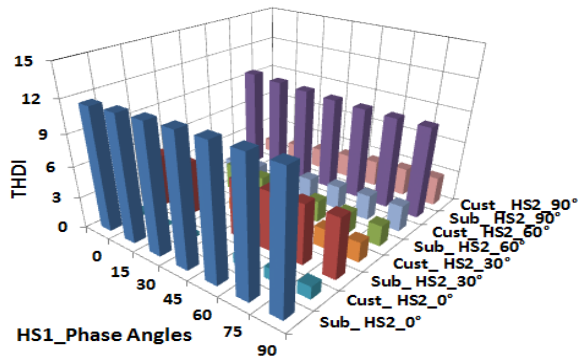

Figure 8. THDI substation vs. customer load point for phase A

## IV. CONCLUSIONS AND REMARKS

This paper investigates a phasor-based method for assessing interactions of multiple harmonic sources in distribution systems. We offer the following observations and conclusions:

1) The new proposed Phasor Harmonic Index, PHI, incorporates more information than the commonly used THDV and THDI indices. PHI considers the phase angles of the distorted voltage and current waveforms in quantifying harmonic content, because phase angle plays a significant role in the interactions between harmonic sources. 2) Phase angles of harmonic sources have complex impacts on the overall harmonic distortion. In some cases, phase angle variations of different harmonic sources result in reduced harmonic impacts. However, phase angle variations that increase harmonic distortion need to be understood. 3) THDV and THDI values are not designed to assess the impact of a single-phase harmonic source on other phases. However, the mutual coupling creates harmonic propagation in all phases. The proposed PHI-I index is helpful in quantifying harmonic distortion in all phases with single-phase harmonic sources present. 4) Harmonic impacts on customer loads and at the substation are evaluated and compared. THD observations show more current distortion at the substation, but more voltage distortion at the customer load. In harmonic studies and in harmonic measurements, harmonic values should be considered throughout the circuit.